\documentclass[12pt,a4paper,reqno]{amsart}
\usepackage{amsfonts,amsthm,latexsym,amsmath,amssymb,amscd,amsmath, epsf}

\usepackage[T2A]{fontenc}
\usepackage[cp1251]{inputenc}
\usepackage{latexsym,amssymb,amsthm}
\usepackage{url,graphics} 
\usepackage[dvips]{graphicx,epsfig} 

\newtheorem{thm}{Theorem}

\newtheorem{prop}[thm]{Proposition}
\newtheorem{lem}[thm]{Lemma}
\newtheorem{defin}{Definition}
\newtheorem{Not}[defin]{Notation}

\newtheorem{rem}{Remark}

\def\q#1.{{\bf #1.}}
\renewcommand\geq{\geqslant}

\newcommand{\be}{\begin{equation}}
\newcommand{\ee}{\end{equation}}

\newcommand{\KK}{\mathbb{K}}

\newcommand {\cal} {\mathcal }
\renewcommand {\hat} {\widehat }
\begin{document}
          \numberwithin{equation}{section}

          \title[$t$-labeled subforests]
					{On commutative algebra associated to $t$-labeled subforests of a graph}

\author[G. Nenashev]
{Gleb Nenashev}
\address{  Department of Mathematics, Stockholm University, SE-106 91 Stockholm, Sweden,}
\email{nenashev@math.su.se}

\begin{abstract}  

For a given graph $G$, we construct an associated commutative algebra, whose dimension is equal to the number of $t$-labeled forests of $G$.

We show that the dimension of the $k$-th graded component of this algebra also has a combinatorial meaning and that its Hilbert polynomial can be expressed through the Tutte polynomial of $G$.
\end{abstract}

\maketitle

\section{Introduction}

The famous matrix-tree theorem of Kirchhoff (see~\cite{Kir} and p.~138 in~\cite{Tut}) claims that the number of spanning trees of a given graph $G$  equals to the determinant of the Laplacian matrix of $G$.
 It is also well known that the number of spanning forests of $G$ or equivalently trees for connected $G$  equals to $T_G(1,1)$ and the number of 
all subforests of $G$ equals to $T_G(2,1)$, where $T_G$ is the Tutte polynomial of $G$ (see e.g.~237 in~\cite{Tut}).

There exist many  generalization of the matrix-tree theorem,  e.g. for directed graphs,  matrix-forest theorems, etc (see~e.g.~\cite{ChK}).
In particular, in~\cite{PSh} A.~Postnikov and B.~Shapiro constructed several algebras associated to~$G$ whose dimensions are
 equal to the number of either spanning trees or forests of~$G$. Below we extend construction of \cite{PSh} to a larger 
class of algebras.

 Given a graph $G$; let as associate commuting variables $\phi_e, e\in G$ to all edges of $G$. For a given positive integer $t\geq 1$, let $\Phi_{G}^{F_t}$ be the algebra generated by $\{ \phi_e: \ e\in G \}$ with relations $\phi_e^{t+1}=0$, for any $e\in G$.

Take any linear order of vertices of $G$. For $i=1,\ldots,n$, Set 
$$X_i=\sum_{e\in G} c_{i,e} \phi_e,$$
 where $c_{i,e}=\pm 1$ for vertices incident to $e$ (for the smaller vertex, $c_{i,e}=1$, for the bigger vertex, is $c_{i,e}=-1$) and $0$ otherwise.
Denote by ${\cal C}_{G}^{F_t}$ the subalgebra of $\Phi_{G}^{F_t}$ generated by $X_1,\ldots ,X_n$.

 Let $\KK$ be some field of characteristic $0$. Consider the ideal $J_{G}^{F_t}$ in the ring $\KK[x_1,\cdots,x_n]$ generated by
$$p_I^{F_t}=\left(\sum_{i\in I} x_i\right)^{tD_I+1},$$
where $I$ ranges over all nonempty subsets of vertices, and $D_I$ is the total number of edges from vertices in $I$ to vertices outside the subset $I$.  Define the algebra   ${\cal B}_{G}^{F_t}$ as the quotient $\KK[x_1,\dots,x_n]/ J_{G}^{F_t}.$

 \begin{Not} 
Fix some  linear order on the edges of $G$. Let $F$ be any a subforest in  $G$.
By $act_G(F)$ denote the number of all externally active edges of $F$, i.e. the number of edges $e\in G\setminus F$ such that subgraph $F+e$ has a cycle and $e$ is the minimal edge in this cycle in the above linear order.

 Denote by $F^+$ the set of edges of the forest $F$ together with externally active edges,
and denote by $F^-=G\setminus F^+$ the set of nonactive edges.
\end{Not}
It is well known that the number of spanning trees and subforests with fixed external activity is independent of the linear order on the set of edges of~$G$. 

For $t=1$ these algebras (denoted by ${\cal B}_{G}^{F}$ and ${\cal C}_{G}^{F}$) were introduced in~\cite{PSh} where the following result was proved.

\begin{thm}[cf.~\cite{PSh}] 
The algebras $B_G^F$ and $C_G^F$ are isomorphic.
Their total dimension as vector spaces over $\KK$ is equal to the number of subforests in the graph $G$.

The dimension of the $k$-th graded component of these algebras  equals
the number of subforests $T$ of $G$ with external activity $|G|-|T|-k$.

\end{thm}

Below we generalize this result for $t>1$, and show that the corresponding dimension coincides with the number of the so-called $t$-labeled trees.
   In Theorem~\ref{tutt} we prove that the Hilbert polynomial of ${\cal B}_{G}^{F_t}$ can be expressed in terms of the Tutte polynomial of $G$. And conversely, in Proposition~\ref{calc} we show that the Tutte polynomial of $G$
can be restored from the Hilbert series of the algebra ${\cal B}_{G}^{F_t}$  for any sufficiently large $t$.

{\bf Acknowledgement.} I am grateful to B~Shapiro for introducing me to this area, and for his editorial help with this text.

\section{$t$-labeled forests}

Consider a finite labelling set containing $t$ diferent labels; each label corresponds to a number from $1$ to $t$. 

\begin{defin}
A spanning forest of the graph $G$ with a label on each edge
is called a $t$-labeled forest.
The weight of a $t$-labeled forest $F$, denoted by $\omega(F)$, is the sum of labels of all its edges. 
\end{defin}

\begin{thm} 
\textsc{(I)} For any graph $G$ and a positive integer $t$, algebras  ${\cal B}_{G}^{F_t}$ and ${\cal C}_{G}^{F_t}$ are isomorphic,
 their total dimension over $\KK$ is equal to the number of $t$-labeled forests in $G$.

\smallskip
\noindent \textsc{(II)} The dimension of the $k$-th graded component of the algebra ${\cal B}_{G}^{F_t}$ is equal to
the number of $t$-labeled forests $F$ of $G$ with the weight $t\cdot (e(G)-act_G(F)) - k$.

\label{tcol}
\begin{proof}
Denote by $\hat{G}$ the graph on $n$ vertices and $t\cdot e(G)$ edges such that each edge of $G$ corresponds to $t$ clones 
in the graph $\hat{G}$, i.e. each edge is substituted by its $t$ copies with labeles $1,2,\ldots,t$. For each edge $e\in G$, its clones $e_1, \dots, e_t\in \hat{G}$  are ordered according to their numbers; 
clones of different edges have the same linear order as the original edges.

Сonsider the following bijection between $t$-labeled forests in $G$ and forests in $\hat{G}$:
each $t$-labeled forest $F\in G$ coresponds to the forest $F'\in \hat{G}$, such that
for each edge $e\in F$, the forest $F'$ has the clone of the edge $e$ whose number is identical with the label of edge $e$ in the forest $F$.
 
Obviously,
$$act_{\hat{G}}(F')=t\cdot act_G(F)+\omega(F)-F,$$
and $e(\hat{G})=t\cdot e(G)$. Since ${\cal B}_{G}^{F_t}$ and ${\cal B}_{\hat{G}}^F$ are the same, the 
Hilbert series of the algebra ${\cal B}_{G}^{F_t}$ coincides with the Hilbert series of the algebra ${\cal B}_{\hat{G}}^F$, 
which settles the second part of Theorem~\ref{tcol}.

To prove the first part of the theorem, observe that ${\cal B}_{G}^{F_t}$ and ${\cal B}_{\hat{G}}^F$ are the same, and algebras 
${\cal C}_{\hat{G}}^{F_t}$ and ${\cal B}_{\hat{G}}^F$ are isomorphic. Thus we must show that algebras  ${\cal C}_{\hat{G}}^{F}$ and ${\cal C}_{G}^{F_t}$ are isomorphic.
 It is indeed true, because for every edge $e\in G$, the elements $\phi_e,\dots ,\phi_e^{t}$ are linearly independent in the algebra $\Phi_{G}^{F_t}$ with coefficients containing no $\phi_e$. 
Also elements $(\phi_{e_1}+\dots+\phi_{e_t}),\dots ,(\phi_{e_1}+\dots+\phi_{e_t})^{t}$  are linearly independent in the algebra $\Phi_{\hat{G}}^{F}$ with coefficients containing no $\phi_{e_1},\dots,\phi_{e_t}$, and  $(\phi_{e_1}+\dots+\phi_{e_t})^{t+1}=0$. Moreover  elements $\phi_{e_i}$ only occur in the sum $(\phi_{e_1}+\dots+\phi_{e_t})$ in the algebra $\Phi_{\hat{G}}^{F}$.
\end{proof}
\end{thm}

\bigskip

 Denote by $c(G)$ the number of connected components of the graph~$G$.
\begin{thm} 
\label{tutt}
Dimension of the $k$-th graded component of ${\cal B}_{G}^{F_t}$ is equal to the coefficient of the monomial $y^{t\cdot e(G) -c(G)+v(G)+1 -k}$ 
in the polynomial $$\left( \frac{y^t-1}{y-1} \right)^{v(G)-c(G)} \cdot T_G \left( \frac{y^{t+1}-1}{y^{t+1}-y}, y^t \right).$$

\begin{proof}

Consider the graph $\hat{G}$ constructed in the proof of Theorem~\ref{tcol}. Set $J_G(x,y):=T_{\hat{G}} (x,y) $; we will use the follow deletion–contraction recurrence for $J_G$,
where  $G-e$ denote  the graph obtained by deleting $e$ from $G$ and $G\cdot e$ is the contraction of $G$ by $e$.

\begin{lem} Polynomial $J_G(x,y)$ satisfies the following:
\begin{enumerate}
	\item If $G$ is empty, then  $J_G(x,y)=1$.
	\item If $e$ is a loop in $G$, then  $J_G(x,y)=y^t J_{G-e}(x,y)$.
	\item If $e$ is a bridge in $G$, then  $J_G(x,y)
	=(y^{t-1} +\dots +1 )\cdot J_{G\cdot e}(x,y) + (x-1)\cdot J_{G-e}(x,y)$.
	\item If $e$ is not a loop or a bridge, then $J_G(x,y)=(y^{t-1} +\dots +1 )\cdot J_{G\cdot e}(x,y) + J_{G-e}(x,y)$.
\end{enumerate}

\label{reclem}
\begin{proof} 

We prove these relations by using the deletion–contraction recurrence for the usual Tutte polynomial
\begin{enumerate}
	\item If $G$ is empty, then  $T_G(x,y)=1$.
	\item If $e$ is a loop in graph $G$, then  $T_G(x,y)=y\cdot T_{G-e}(x,y)$.
	\item If $e$ is a bridge in graph $G$, then  $T_G(x,y)=x\cdot T_{G-e}(x,y)$.
	\item If $e$ is not a loop or a bridge, then $T_G(x,y)=T_{G\cdot e}(x,y) + T_{G-e}(x,y)$.
\end{enumerate}

\noindent {\bf 1.} Graph ${\hat{G}}$ is also empty, hence   $J_G(x,y)=T_{\hat{G}} (x,y) =1$.
\smallskip

\noindent {\bf 2.} Clones $e_1,\dots , e_t$ is also loops in graph ${\hat{G}}$, therefore 
$T_{\hat{G}} (x,y) =y^t\cdot T_{\hat{G}- \{ e_1,\dots , e_t\}} (x,y) = y^t\cdot T_{\hat{G- e}}(x,y)$, 
hence   $J_G(x,y)=y^t\cdot T_{\hat{G- e}}(x,y) =J_{G-e}(x,y)$. 
\smallskip

\noindent {\bf 3.} We calculate our polynomial using the deletion–contraction recurrence for the Tutte polynomial.
$$J_G(x,y)=T_{\hat{G}} (x,y)=T_{\hat{G}\cdot e_k} (x,y) + T_{\hat{G} - e_t} (x,y) =$$
$$y^{t-1}\cdot T_{\hat{G\cdot e}}(x,y) + T_{\hat{G} - e_t} (x,y) = $$
$$y^{t-1}\cdot T_{\hat{G\cdot e}}(x,y) + T_{(\hat{G} - e_t)\cdot e_{t-1}} (x,y) + T_{\hat{G} - e_t- e_{t-1}} (x,y)=$$
$$(y^{t-1}+y^{t-2})\cdot T_{\hat{G\cdot e}}(x,y) + T_{\hat{G} - e_t- e_{t-1}} (x,y)=$$
$$\cdots $$
$$(y^{t-1}+y^{t-2}+\ldots +y)\cdot T_{\hat{G\cdot e}}(x,y) + T_{\hat{G} - e_t- \ldots -e_{2}} (x,y)=$$
$$(y^{t-1}+y^{t-2}+\ldots +y)\cdot T_{\hat{G\cdot e}}(x,y) + x\cdot T_{\hat{G} - e_t- \ldots -e_{1}} (x,y)=$$
$$(y^{t-1}+y^{t-2}+\ldots +y)\cdot T_{\hat{G\cdot e}}(x,y) + x\cdot T_{\hat{G- e}} (x,y)=$$
$$(y^{t-1}+y^{t-2}+\ldots +y+1)\cdot T_{\hat{G\cdot e}}(x,y) + (x-1)\cdot T_{\hat{G- e}} (x,y)=$$
$$(y^{t-1} +\dots +1 )\cdot J_{G\cdot e}(x,y) + (x-1)\cdot J_{G-e}(x,y).$$
\smallskip
\noindent {\bf 4.} It is similar to { 3}, but now we have 
$$T_{\hat{G} - e_k- \ldots -e_{2}} (x,y)= T_{(\hat{G} - e_k- \ldots -e_{2})\cdot e_1} (x,y) +T_{\hat{G} - e_k- \ldots -e_{2}-e_{1}} (x,y)=$$ $$T_{\hat{G\cdot e}}(x,y) +  T_{\hat{G- e}} (x,y),$$
since in this case edge $e_1$ is not a bridge in $\hat{G} - e_k- \ldots -e_{2}$.

\end{proof}
\end{lem}

Now let us rewrite $J_G(x , y)$ in terms of $t$-labeled forests using the deletion–contraction recurrence for $J_G(x , y)$ in the above fixed linear order of edges of $G$.  Obviously,
 the edges by which we contract the graph constitute a forest. Therefore,  $J_G(x , y)=\sum_{F_u} a({F_u})$, where $a({F_u})$ depends only   
on $G$ and the forest $F_u$. Now rewrite the latter equality in terms of $t$-labeled forests. When we contract edge $e$ in $G$, the term $y^{k-1}$ in the factor $(y^{t-1}+y^{t-2}+\ldots +y+1)$ corresponds to the choice  the $k$-th label for edge $e$, i.e. we have $J_G(x , y)=\sum_{F}y^{\omega(F)-|F|} b({F})$.
It remains to calculate $b(F)$. An edge for $t$-labeled forest $F$ is a loop if and only if it is active, and the number of edges 
which are bridges in our recursion equals to $c(F)-c(G)=(v(G)-1-|F|)-c(G)= (v(G)-c(G)) -1-|F|$. Therefore, we have

$$J_G(x , y)=\sum_{F}y^{\omega(F)-|F|}\cdot y^{t\cdot act_G(F)} \cdot (x-1)^{(v(G)-c(G))-1-|F|},$$ 
$$J_G(x , y)=\sum_{F}y^{\omega(F)-|F|+ act_G(F)} \cdot (x-1)^{(c(G)-v(G))-1-|F|},$$
$$J_G(1+\frac{1}{y} , y)=\sum_{F}y^{\omega(F)-|F|+ t \cdot act_G(F)} \cdot (\frac{1}{y})^{(c(G)-v(G))-1-|F|},$$
$$J_G(1+\frac{1}{y} , y)=\sum_{F}y^{\omega(F)-|F|+ t \cdot act_G(F)-((c(G)-v(G))+1+|F|},$$
$$J_G(1+\frac{1}{y} , y)=\sum_{F}y^{\omega(F)+ t \cdot act_G(F)-c(G)+v(G)+1}.\eqno (*)$$
By Theorem~\ref{tcol} the dimension of the $k$-th graded component of algebra ${\cal B}_{G}^{F_t}$ equals the number of $t$-labeled forests $F$ of $G$ with weight $t\cdot (e(G)-act_G(F)) - k$. Then the dimension of the $k$-th graded component is equal to the coefficient of the monomial $y^{t\cdot e(G) -k -c(G)+v(G)+1}$ in  polynomial $J_G(1+\frac{1}{y} , y)$.

\begin{lem}
\label{=}
$$J_G(1+\frac{1}{y} , y)= \left( \frac{y^t-1}{y-1} \right)^{v-c(G)}\cdot  T_G \left( \frac{y^{t+1}-1}{y^{t+1}-y}, y^t \right). $$
\begin{proof}
Conditions $1,2$ and $4$ of Lemma~\ref{reclem} hold for polynomial $\left( \frac{y^t-1}{y-1} \right)^{v-c(G)}\cdot  T_G \left( \frac{y^{t+1}-1}{y^{t+1}-y}, y^t \right)$. 
Now we can check the  $3$-rd condition. Set $$z:=(y^{t-1}+\ldots+1)=\frac{y^t-1}{y-1},$$ then, $\frac{y^{t+1}-1}{y^{t+1}-y} = \frac{1}{zy}+1.$ We have 
$$z^{v(G)-c(G)}\cdot  T_G \left( \frac{1}{zy}+1, y^t \right)=$$
$$z^{v(G)-c(G)}\cdot  T_{G\cdot e} \left( \frac{1}{zy}+1, y^t \right) + \frac{1}{zy}\cdot z^{v(G)-c(G)}\cdot  T_{G-e} \left( \frac{1}{zy}+1, y^t \right)=$$
$$(y^{t-1}+\ldots+1)\cdot z^{v(G\cdot e)-c(G\cdot e)}\cdot  T_{G\cdot e} \left(\frac{1}{zy}+1, y^t \right) +$$ $$ \frac{1}{y}\cdot z^{v(G-e)-c(G-e)}\cdot  T_{G-e} \left( \frac{1}{zy}+1, y^t \right).$$

Hence, the  $3$-rd condition holds as well.
Therefore, if we calculate these polynomials using the recursion method we get the same results, hence, these polynomials coincide.
\end{proof}
\end{lem}

This settles Theorem~\ref{tutt}.
\end{proof}
\end{thm}

\begin{prop}
\label{calc}
For any positive integer $t\geq n$, it is possible to restore the Tutte polynomial of any connected graph $G$ on $n$ vertices 
knowing only the dimensions of each graded component of the algebra ${\cal B}_{G}^{F_{t}}$.

\begin{proof}
Choose a integer $t\geq n$.
By Theorem~\ref{tcol} we know that the degree  of  the maximal non empty graded  component of ${\cal B}_{G}^{F_{t_n}}$
equals to the maximum of $t\cdot(e(G)-act_G(F))-\omega(F)$ taken over $F$. It attains its maximal value for the empty forest (i.e. $F=\emptyset$). 
Then we know the value of $t\cdot e(G)$, hence, we know the number of edges of the graph $G$.

By Theorem~\ref{tutt} we also know the polynomial $$\left( \frac{y^t-1}{y-1} \right)^{v(G)-c(G)} \cdot T_G \left( \frac{y^{t+1}-1}{y^{t+1}-y}, y^t \right),$$ because $G$ is connected (i.e. $c(G)=1$). This polynomial equals to $$\sum_{F}y^{\omega(F)+ t \cdot act_G(F)-c(G)+v(G)+1}=\sum_{F}y^{\omega(F)+ t \cdot act_G(F)+v(G)},$$  where the summation is taken over all $t$-labeled forests (see eq.~(*) and Lemma~\ref{=}).
Rewriting it in terms of the usual subforests,  we can calculate 
$$\sum_{F_u}(y+\ldots+y^t)^{|F_u|}\cdot y^{ t \cdot act_G(F_u)+v(G)},$$
Hence, we also know the sum
$$\sum_{F_u}(1+\ldots+y^{(t-1)})^{|F_u|}\cdot y^{ |F_u|+t \cdot act_G(F_u)}.\eqno (**)$$

Since $|F_u|<t$, then we can compute the number of usual subforests with a fixed pair of parameters $|F_u|$ and $act_G(F_u)$.
Consider the monomial of minimal degree in polynomial~$(\oplus)$, and present it in the form $s\cdot y^m$. Observe that $s$ is the number of subforests $F_u$ s.t. $F_u\equiv m\ (mod\ t)$ and with $act_G(F_u)=\left[\frac{m}{t}\right].$ Remove from the polynomial~$(**)$ all summands for these subforests, and repeat this operation until we get $0$.

It is well known that
$T_G(x,y)=\sum_{a,b} \#\{ F_u : |F_u|=a, \ act(F_u)=b\}\cdot (x-1)^{n-1-a} \cdot y^b$. Therefore since we know the number of usual subforeests with any fixed number of edges and any fixed extrenal
activity, we know the whole Tutte polynomial.

\end{proof}
\end{prop}

\begin{rem}
It is possible to obtain similar results for $t$-labeled trees except for Proposition~\ref{calc}. But in our opinion such results are not very interesting, because the number of edges in every tree is the same.
\end{rem}


\end{document}